\documentclass{article}
\usepackage{amsmath}
\usepackage{amssymb}
\usepackage{amsfonts}
\usepackage{graphicx}
\usepackage{multicol}
\usepackage{blindtext}
\usepackage{enumitem}
\usepackage[all]{xy}

\newtheorem{theorem}{Theorem}[section]
\newtheorem{proposition}[theorem]{Proposition}
\newtheorem{lemma}[theorem]{Lemma}

\newtheorem{remark}[theorem]{Remark}

\makeatletter

\@addtoreset{equation}{section}
\makeatother

\newcommand{\ve}{\varepsilon}

\newcommand{\R}{\ensuremath{\mathbb{R}}}

\def\com{\mathbb C}
\def\r{\mathbb R}

\title{A Finsler metric of constant Gauss curvature $K = 1$ on 2-sphere \footnote{to appear in An. \c Stiin\c t. Univ. Al. I. Cuza Ia\c si Mat. (N.S.)}}

\author{By  I. Masca, S. V. Sabau, H. Shimada}
\date{}
\begin{document}
\maketitle

\begin{abstract}
We construct a concrete example of constant Gauss curvature $K = 1$ on the 2-sphere having all geodesics closed and of same length.\\


{\it Subclass}: 53C60, 53C22
\end{abstract}

\section{Introduction}
Zoll surfaces are Riemannian metrics all of whose geodesics are simple closed curves of equal length (see for example
\cite{5}, \cite{6} and many other authors). They are natural generalizations of the metric on the round  sphere, having the Gauss curvature of Zoll metrics not constant.

More precisely, a Zoll surface of rotation on $(S^2, g)$ is a Riemannian metric on the 2-sphere with canonical spherical coordinates $(r, \theta) \in \lbrack0, \pi\rbrack \times (0, 2\pi\rbrack$ given by 
\begin{equation}\label{Zoll metric }
g=[1+h(\cos r)]^2 dr\otimes dr+\sin^2\ r d\theta\otimes d\theta,
\end{equation}
where $h:[-1,1]\to (-1,1)$ is a smooth function such that 
\begin{enumerate}[label=(\roman*)]
 \item 
 $h$ is odd function, i.e. $h(-x) = -h(x)$, for all $x \in [-1,1]$,
 \item
 $h(-1) = h(1) = 0$.
\end{enumerate}
Observe that $h(0) = 0$ from condition (i).

It can be checked that this Riemannian has all geodesics closed, of same length, and with Gauss curvature 
\begin{equation}
 G(r) = \frac{1 + h(\cos r) - \cos r\cdot h'(\cos r)}{[1 + h(\cos r)]^3}
\end{equation}
(see \cite{5} or \cite{1}  for detailed computations).

Zoll metrics have many remarkable geometrical properties (see \cite{5}). We recall here only that the manifold of oriented geodesics, say $M$, of a Zoll surface $(S^2, g)$ is a smooth manifold diffeomorphic to $S^2$.\\

One problem to ask is {\it what kind of natural geometrical structures are carried by the manifold of geodesics $M$ of a Zoll surface $(S^2, g)$?}\\

Some of the answers are already known:
\begin{enumerate}[label=(\roman*)]
 \item 
 $M \simeq S^2$ can be endowed with a symplectic structure, or it can be regarded as a Riemannian manifold with an averaged Riemannian metric (see \cite{5}) for details).
 \item 
 $M$ also carries a Finsler structure of constant Gauss curvature $K = 1$ (see \cite{3}, \cite{1} for details).
\end{enumerate}

The existance of a such Finsler structure on the manifold of geodesics was shown for the fist time by Bryant in \cite{3} by using exterior differential systems. He also constructed a Finsler metric of constant Gauss curvature $K = 1$ on $S^2$ with the supplementary condition of projectively flatness by using an approach borrowed from algebraic geometry. However, this metrics (called today Bryant metrics) do not come from Zoll metrics (as pointed out by Bryant himself in \cite{3}, p. 254, bottom).

The concrete construction of the Finsler surfaces of constant Gauss curvature $K = 1$ on $S^2$ using positively curved Zoll metrics, was studied very recently 
in \cite{1}, where one can find a lot of information about these Finsler metrics on $S^2$ as well as the parametric equation of the indicatrix. However the concrete form of the fundamental function $F$ was not yet explicitely obtained.

Without being trivial, observe that the Zoll metrics constructed using the following polynomial function:
\begin{eqnarray}\label{ex1 h}  h(x) =
\begin{cases}
 \varepsilon x(1 - x^2), \qquad 0 < \varepsilon < \frac{1}{2}\\
 x(1 - x^2)^n, \qquad n > 1
\end{cases}
\end{eqnarray}
has the property $G(x) > 0$, for $x \in [-1, 1]$, hence each of these functions induces a $K = 1$ Finsler metric on the manifold of geodesics $M$.

The aim of the present paper is to {\it find the explicit form of the fundamental function $F$ of such a Finsler metric} in the simplest case of 
\begin{equation}\label{(1)}
 h: [-1, 1] \rightarrow \r, \qquad h(x) = \varepsilon x(1 - x^2), \qquad 0 < \varepsilon < \frac{1}{2}.
\end{equation}

We will study the solutions of this equation, find the precise positive real solution that corresponds to the Finsler metric and write it explicitely in a formula with radicals.

It will be seen that this Finsler metric is actually an $(\alpha, \beta)$-metric, albeit a quite complicated one.

This Finsler metric, even quite complicated when explicitely written in coordinates has some remarkable properties
\begin{enumerate}[label=(\roman*)]
 \item 
 it is of constant Gauss curvature $K = 1$ on $S^2$,
 \item
 it has all geodesics closed and of same length,
 \item
 it is a Finsler surface of revolution, i.e. it is rotationally invariant,
 \item
 the unit speed geodesics of this Finsler structure, 
emanating from the some point $p \in M$ intersect at the same distance from $p$ in the some point, i.e. the cut locus of a point $p \in M$ is a point.
\end{enumerate}

In the Riemannian case, a metric with these properties  is unique and it must be the canonical metric of the round sphere, but in the Finsler case, we have a family of such Finsler metrics, each of them corresponding to an odd function $h: [-1, 1] \rightarrow (-1, 1)$ with some supplementary properties such as $G(r) > 0$. We single out one such Finsler metric in the present paper. 

We point out that in the case of the function \eqref{ex1 h} with $n = 2$, the implicit indicatrix equation is an 8-th order algebraic equation in $F$. The famous Abel-Ruffini theorem says that there is no solution in radicals to a general equation of degree 5 or higther, with arbitrary coefficients, hence a similar approch with the present paper is impossible.
New algebraic and analytical notices are needed and we will leave this case for the future research.

\section{The settings and notations}	
	
Let us start with
the odd function $h$ in \eqref{(1)}, and observe that the Gauss curvature of the corresponding Zoll metric $g$ on $S^2$ is 
\begin{equation}\label{G for g}
 G(x) = -\frac{2\varepsilon x^3 + 1}{(\varepsilon x^3 - \varepsilon -1)^3},
\end{equation}
that is obviously positive on $[-1, 1]$, for $0 < \varepsilon < \frac{1}{2}$.\\

Moreover, let us recall the following essential result from \cite{1}.

	\begin{theorem}\label{thm: indicatrix param eqs}
			The parametric equations of the indicatrix of the corresponding Finsler metric of constant flag curvature $K=1$ are 
		\begin{equation}\label{param eq indic ver1}
		\begin{split}
		 v_1(r) & = \pm\frac{1}{\cos R}\cdot \sqrt{\sin^2r-c^2}\\
		 v_2(r) & = \frac{1+h(\cos r)}{\cos r}\\	
		& -\sqrt{\sin^2r-c^2}\int_{r_c}^r \frac{\sin s}{\cos^2s}\Bigl[1- 
		\frac{\cos s\cdot h'(\cos s)}{1+h(\cos s)} \Bigr]
		\Bigl[ \frac{1+h(\cos s)}{\sqrt{\sin^2s-c^2}} \Bigr] ds.
		\end{split}
		\end{equation} 	
	\end{theorem}

	Observe that the sign of $v_2$ is different from \cite{1}, but this makes no difference for our research. We will compute now the algebraic equation for the Finsler metric $F$. This result follows from a more general result in \cite{1}, but we compute it here directly for the sake of completeness. 

        Let us start by remarking that, by using the function $h$ in \eqref{(1)}, and some elementary straightforward computations, the second equation in \eqref{param eq indic ver1} reads
        \begin{equation}
          \begin{split}
          v_2 & =\frac{1+\varepsilon \sin^2r \cos r}{\cos r}\left[
\frac{\sin^2 r-c^2}{\cos^2 R\cos r}+2\varepsilon(\sin^2r-c^2)
\right]\\
&= \frac{\cos r}{\cos^2R}-\varepsilon(\sin^2r-c^2)+\varepsilon c^2,
\end{split}
          \end{equation}
          where $c=\sin R$. This leads to
          \begin{equation}
\frac{\cos r}{\cos ^2 R}=v_2 +\varepsilon(\sin^2r-c^2)-\varepsilon c^2,
            \end{equation}
and from here, it follows that for this choice of $h$, the implicit equation of the indicatrix curve of $F$ is 
\begin{equation}\label{ex1:implicit formula}
\frac{1-v_1^2}{\cos^2R}=(v_2+\ve v_1^2\cos^2R  -\ve c^2)^2,
\end{equation}
where $(R, \theta) \in \big[-\frac{\pi}{2}, \frac{\pi}{2}\big]\times [0, 2\pi)$ are the coordinates on the manifold of geodesics $M \simeq S^2$, $(v_1, v_2)$ in the coordinates of the fiber $T_{(R, \theta)}M$, and $c = \sin R$.\\

In order to obtain the algebraic equation for the fundamental function of this Finsler space $F$, we simply substitute $v_i$ by $\dfrac{v_i}{F}$, $i \in\{ 1, 2\}$, since the indicatrix equation is given by $F=1$. It follows 
\begin{equation}
\frac{F^2-v_1^2}{\cos^2 R}=\frac{(\ve c{^2F^2-\ve v_1^2\cos^2R-Fv_2})^2}{F^2},
\end{equation}
and after some computations, and the use of $c = \sin R$, we put this algebraic equation in the canonical form 
\begin{equation}\label{eq1 in F simple}
A F^4+BF^3+C F^2+DF+E=0,
\end{equation}
where 
\begin{equation}\label{coeficienti1}
\begin{split}
  A & :=1-\ve^2(1-c^2)c^4 
  \\
  B & :=2\ve c^2(1-c^2)v_2 
  \\
C & :=(2c^6\ve^2-4c^4\ve^2+2c^2\ve^2-1)v_1^2-(1-c^2)v_2^2\\
D & :=-2\ve(1-c^2)^2 v_1^2v_2 
\\
E & :=-\ve^2v_1^4(1-c^2)^3. 
\end{split}
\end{equation}
We need to find a real positive solution of this 4-th order algebraic equation. Obviously, by using the substitution $c=\sin R$, there formulas can be written in terms of the coordinate $R$, but we do not  need the explicit form for these.

A glance at the  coefficients of the equation \eqref{eq1 in F simple} gives the following.

\begin{lemma}\label{semnele}
  For any $\varepsilon \in (0, \frac{1}{2})$,  the coefficients $A$, $C$, $E$ of \eqref{eq1 in F simple} satisfy
  \begin{equation}
A>0,\quad C<0,\quad E<0.
    \end{equation}

\end{lemma}

If we think the equation \eqref{eq1 in F simple} as defininig a Minkowski type metric in each tangent plane $T_{(R,\Theta)}M$, then it is clear that its coefficients are constants and hence we can  consider the function $f: \R \longrightarrow \R$, $f(x) = Ax^4 + Bx^3 + Cx^2 + Dx + E$. Observe that
\begin{enumerate}
	\item $f(0)=E<0$, and 
	\item $\lim_{x\to \pm \infty}f(x)=+\infty$. 
\end{enumerate}
Since $f$ is a smooth (hence continuous) function, we obtain
\begin{proposition}
	The algebraic equation 
	\eqref{eq1 in F simple} has at least two real solutions, one non-positive and one non-negative. 
\end{proposition}

Taking into account now the Vieta's formulae for \eqref{eq1 in F simple}, that is
\begin{equation}\label{Vieta 1}
  F_1F_2F_3F_4=\dfrac{E}{A}<0,
  \end{equation}
it
 results that there are only two situations possible:

\begin{enumerate}
	\item there are two real solutions of opposite signs $(+,-)$, and two complex conjugate solutions, or 
	\item all four solutions are real, having the  signs $(+---)$ or $(+++-)$.

\end{enumerate}

\section{The solution}
	
Let us recall that there are several equivalent methods available to solve 4-th order algebraic equation, like Ferrari's method, Galois' method and maybe others.

We will work in the following with the Galois' method (see \cite{8} for a very simple  exposition). It is well-known that over $\com$, the algebraic equation \eqref{eq1 in F simple} has four solutions $F_k$, $k \in\{1,2, 3, 4\}$.

We start by making the substitution 
\begin{equation}\label{F to X}
F=X-\frac{B}{4A},
\end{equation}
such that the equation \eqref{eq1 in F simple} reduces to
\begin{equation}\label{eq1 in X simple}
X^4+\alpha X^2+\beta X+\gamma=0,
\end{equation}
where $\alpha$, $\beta$, $\gamma$ are some constants that can be written in terms of the initial coefficients of \eqref{eq1 in F simple}, but we don't need their explicit form.

This is called the {\it reduced quartic polynomial equation}. It is clear that over $\com$, this equation must have four solutions.

\begin{remark}
We mention that some more computation shows that the coefficients of the reduced equation \eqref{eq1 in X simple} are polynomials in the fiber coordinates $(v_1,v_2)$ with the coefficients depending only on $R$, that is constants in a tangent plane $T_{(R,\Theta)}M$. More precisely, we get  
\begin{equation}\label{alpha,beta,gamma}
\begin{split}
& \alpha=\alpha_1v_1^2+\alpha_2v_2^2\\
& \beta=\beta_{12}v_1^2v_2+\beta_2v_2^3\\
&\gamma=\gamma_1v_1^4+\gamma_2v_2^4+\gamma_{12}v_1^2v_2^2.
\end{split}
\end{equation}
\end{remark}

Let us define the following expressions
\begin{equation}
\begin{split}
z_1 & := (X_1+X_2)(X_3+X_4)\\
z_2 & := (X_1+X_3)(X_2+X_4)\\
z_3 & := (X_1+X_4)(X_2+X_3),
\end{split}
\end{equation}
where $X_k$, $k\in\{1,2,3,4\}$, are the solutions (some of them might be complex) of the quartic equation \eqref{eq1 in X simple},
and let us consider the cubic equation having these expressions as solutions, that is
\begin{equation}
  (z-z_1)(z-z_2)(z-z_3)=0.
\end{equation}

Some computation shows this can be written as 
\begin{equation}\label{Galois resolvent}
  z^3-2\alpha z^2+(\alpha^2-4\gamma)z+\beta^2=0,
\end{equation}
where $\alpha$, $\beta$, $\gamma$ are the coefficients of \eqref{eq1 in X simple}.

This equation is called the {\it resolvent equation} of the equation \eqref{eq1 in X simple}. 

Applying Cardano's formulas (see \cite{8}) for the resolvent equation \eqref{Galois resolvent}, we get the solutions
\begin{eqnarray}\label{resolv sols 1}
z_1 & = & {P} + {Q} + \frac{2\alpha}{3} \\\nonumber
  z_2 & = & \omega{P} + \omega^2{Q} + \frac{2\alpha}{3} = -\frac{1}{2}\Big({P} + {Q} - \frac{4\alpha}{3}\Big) + \frac{i\sqrt{3}}{2}({P} - {Q})
  \\\nonumber
z_3 & = & \omega^2{P} + \omega{Q} +\frac{2\alpha}{3} = -\frac{1}{2}\Big({P} + {Q} - \frac{4\alpha}{3}\Big) - \frac{i\sqrt{3}}{2}({P} - {Q}).
\end{eqnarray}
Here, $\omega$ is the cubic root of the unity, and
\begin{eqnarray}\label{P-uri}
{P} & = & \sqrt[3]{\mathcal A_{(6)} + \sqrt{\mathcal B_{(12)}}},\\\nonumber
{Q} & = & \sqrt[3]{\mathcal A_{(6)} - \sqrt{\mathcal B_{(12)}}},
\end{eqnarray}
where $\mathcal A_{(6)}$ and $\mathcal B_{(12)}$ are some homogeneous polynomials of orders 6 and 12 in the fiber coordinates $(v_1,v_2)$, respectively. 

The resolvent is used in the following way to obtain the solutions of the original quartic equation.

If we consider 
\begin{equation}
0=(X_1+X_2)+(X_3+X_4)\quad \textrm{and}\quad z_1=(X_1+X_2)(X_3+X_4)
\end{equation}
it follows 
\begin{equation}
\begin{cases}
& X_1+X_2=\sqrt{-z_1}\\
& X_3+X_4=-\sqrt{-z_1}.
\end{cases}
\end{equation}

Similarly we get 
\begin{equation}
\begin{cases}
& X_1+X_3=\sqrt{-z_2}\\
& X_2+X_4=-\sqrt{-z_2},
\end{cases}
\end{equation}
and 
\begin{equation}
\begin{cases}
& X_1+X_4=\sqrt{-z_3}\\
& X_2+X_3=-\sqrt{-z_3}.
\end{cases}
\end{equation}

By putting all these equations together, we obtain a linear system of six equations with four unknown that can be studied by elementary means. Indeed, the coefficients matrix is 
\begin{equation}
\begin{pmatrix}

1 & 1 & 0 & 0 \\
0 & 0 & 1 & 1 \\
1 & 0 & 1 & 0 \\

0 & 1 & 0 & 1 \\
1 & 0 & 0 & 1 \\
0 & 1 & 1 & 0 
\end{pmatrix}
\end{equation}
and its rank is clearly 4. It is elementary to see that the solutions can be obtained by solving the linear system 
\begin{equation}
\begin{pmatrix}
1 & 1 & 0 & 0 \\
1 & 0 & 1 & 0 \\
1 & 0 & 0 & 1 \\
0 & 1 & 1 & 0 
\end{pmatrix}
\begin{pmatrix}
X_1\\
X_2\\
X_3\\
X_4
\end{pmatrix}=
\begin{pmatrix}
\sqrt{-z_1}\\
\sqrt{-z_2}\\
\sqrt{-z_3}\\
-\sqrt{-z_3}
\end{pmatrix}.
\end{equation}

The solution of this system reads
\begin{equation}\label{X by z}
  \begin{split}
X_1 & =  \frac{1}{2}\left[\sqrt{-z_1}+\sqrt{-z_2}+\sqrt{-z_3}\right] \\
X_2 & =  \frac{1}{2}\left[\sqrt{-z_1}-\sqrt{-z_2}-\sqrt{-z_3} \right] \\
X_3 & =  \frac{1}{2}\left[- \sqrt{-z_1}+\sqrt{-z_2}-\sqrt{-z_3} \right]\\
X_4 & =  \frac{1}{2}\left[-\sqrt{-z_1}-\sqrt{-z_2}+\sqrt{-z_3} \right]. 
\end{split}
\end{equation}

We obtain the main result.
\begin{theorem}
 The fundamental function $F$ of the $K = 1$ Finsler metric on $S^2$ induced by $h$ in \eqref{(1)} is 
 \begin{equation}\label{Gauss1}
   \begin{split}
     F & = \frac{1}{2}\Bigg[\sqrt{-\left({P} + {Q} + \frac{2\alpha}{3}\right)}\\
     &
     +
     \sqrt{-\left({P} + {Q} - \frac{4\alpha}{3}\right) + \sqrt{\left({P} + {Q}-\frac{4\alpha}{3}\right)^2+3(P-Q)^2}}
    \  \Bigg]
     - \frac{B}{4A},
   \end{split}
\end{equation}
where $P$, $Q$ are given in the formulas \eqref{P-uri},
$\alpha$ in \eqref{alpha,beta,gamma}, and 
$A$, $B$ in \eqref{coeficienti1}, respectively. 
\end{theorem}
Indeed,
  observe that the Vieta's formulas for the resolvent equation \eqref{Galois resolvent} give
  $$
z_1z_2z_3=-\beta^2<0,
$$
that is, if all the solutions of the resolvent \eqref{Galois resolvent} are real, the following cases are possible
\begin{enumerate}
\item if there is only a real solution, this solution is negative, or
  \item if all three solutions are real, then either all of them are negative, or
  two solutions are positive and one negative. 
  \end{enumerate}

  However, the case when two solutions are positive and one negative cannot happen because of the following argument. Let us assume, for instance, that
  $$
z_1>0,\quad z_2>0,\quad z_3<0.
  $$
  Then,
  \begin{enumerate}
\item if all solutions are different, i.e. $z_1\neq z_2$, then from \eqref{X by z} it follows that all solutions $X_i$, $i\in \{1,2,3,4\}$ are complex, i.e. the equation \eqref{eq1 in X simple} cannot have real solutions, and hence the same for the  equation \eqref{eq1 in F simple}. Of course this is not possible since we know from geometrical reasons (see \cite{1}) that there exists a real Finsler metric of constant flag curvature $K=1$ whose indicatrix is given by \eqref{param eq indic ver1} for the function $h$ chosen by us, or
\item there is a double solution $z_1=z_2$. Then \eqref{X by z} implies that $X_1$, $X_4$ are complex conjugate solutions and $X_2=X_3$ is a double real solution that must be positive since it is given by radicals. Using \eqref{F to X} it follows that same thing can be said about $F_i=X_i-\dfrac{B}{4A}$, $i\in\{1,2,3,4\}$, as well, namely,  $F_1$, $F_4$ are complex conjugate solutions and $F_2=F_3$ is a double real solution. 

  However, this implies  $F_1F_2F_3F_4>0$ that is contradiction to \eqref{Vieta 1}.
    \end{enumerate}

We conclude that only possible cases are 
\begin{enumerate}
\item if there is only a real solution, this solution is negative, or
  \item if all three solutions are real, then all of them are negative.
  \end{enumerate}

Taking now into account formulas \eqref{X by z}, it follows that 
\begin{equation}
   \begin{split}
     X_1 & = \frac{1}{2}\Bigg[\sqrt{-\left({P} + {Q} + \frac{2\alpha}{3}\right)}\\
     &
     +
     \sqrt{-\left({P} + {Q} - \frac{4\alpha}{3}\right) + \sqrt{\left({P} + {Q}-\frac{4\alpha}{3}\right)^2+3(P-Q)^2}}
    \  \Bigg]
       \end{split}
\end{equation}
is the solution of \eqref{eq1 in X simple} that we are looking for. 

By using the transformation \eqref{F to X}, the conclusion follows. 

\begin{remark}
It worth pointing out that the Finsler metric $F$ in \eqref{Gauss1} belongs to the family of so-called $(\alpha, \beta)$-metrics.
\end{remark}
Indeed, if we consider $(\alpha, \beta)$ such that $v_1 := \alpha^2 - \beta\cos R$ and $v_2 := \beta$, after some computations we get that $\mathcal{A}_{(6)} $ and $\mathcal{B}_{(12)} $ are homogeneous polynomials of degree six and twelve in $\alpha$ and $\beta$, hence the Finsler metric in \eqref{Gauss1} can be written as an $(\alpha,\beta)$ metric.
This fact can be regarded as a consequence of a more general result in \cite{7}. Please do not confuse the  $\alpha$ and $\beta$ used in this Remark with the $\alpha$ and $\beta$ in \eqref{alpha,beta,gamma}.

\section{Why is this Finsler metric of constant flag curvature?}

One legitimate question to ask if the Finsler metric we have computed in the present paper is indeed of constant flag curvature $K=1$.

Let us recall that a  Finsler metric on a real smooth, $n$-dimensional manifold
$M$ is a function $F:TM\to \left[0,\infty \right)$ that is positive and
smooth on $\widetilde{TM}=TM\backslash\{0\}$, has the {\it homogeneity property}
$F(x,\lambda v)=\lambda F(x,v)$, for all $\lambda > 0$ and all 
$v\in T_xM$, having also the {\it strong convexity} property that the
Hessian matrix
\begin{equation}\label{hess F}
g_{ij}=\frac{1}{2}\frac{\partial^2 F^2}{\partial y^i\partial y^j}
\end{equation}
is positive definite at any point $u=(x^i,y^i)\in \widetilde{TM}$.

The fundamental function $F$ of a Finsler structure $(M,F)$ determines and it is determined by the (tangent) {\it indicatrix}, or the total space of the unit tangent bundle of $F$, namely
\begin{equation*}
\Sigma_F:=\{u\in TM:F(u)=1\}
\end{equation*}
which is a smooth hypersurface of $TM$. At each $x\in M$ we also have the {\it indicatrix at x}
\begin{equation*}
\Sigma_x:=\{v\in T_xM \ |\  F(x,v)=1\}=\Sigma_F\cap T_xM
\end{equation*}
which is a smooth, closed, strictly convex hypersurface in
$T_xM$. 

To give a Finsler structure $(M,F)$ is therefore equivalent to giving a smooth
hypersurface $\Sigma\subset TM$ for which the canonical projection
$\pi:\Sigma\to M$ is a surjective submersion and having the property
that for each $x\in M$, the $\pi$-fiber $\Sigma_x=\pi^{-1}(x)$ is
strictly convex including the origin $O_x\in T_xM$.

In the present paper we have started with a Zoll metric $(S^2,g)$ given in \eqref{Zoll metric }
determined by the odd function $h$ in \eqref{(1)}. It is easy to see that the Riemannian metric $g$ given in \eqref{Zoll metric } with $h$ in \eqref{(1)} is indeed a Zoll metric. This follows immediately from \cite{5}. Indeed, it is elementary to see that the Darboux theorem 4.11 in \cite{5} is satisfied by the function $f=1+h$.  

The fundamental geometrical property of Zoll metrics is that they have all geodesics closed and of same length. This implies that the space of geodesics, that is the topological space where a point represents an oriented geodesic, is a smooth differentiable manifold. Observe that this property is not true in general for an arbitrary Riemannian metric. Therefore it is natural to construct a Finsler metric on the manifold of geodesics $M$ of a Zoll metric $(S^2,g)$ such that the indicatrices of this Finsler metric coincide with the closed geodesics of $g$. 

If we denote by $U^gS^2=\Sigma$ the unit sphere bundle of the Zoll metric $(S^2, g)$, then on $\Sigma$ there exists a canonical $g$-orthogonal coframe $\{\theta^1,\theta^2,\theta^3\}$ uniquelly determined by the Riemannian metric $g$. We recall here this construction for the sake of completeness. Let us consider a $g$-orthonormal basis $f_1,f_2$ of the tangent space $T_xS^2$ and denote by $\alpha^1$, $\alpha^2$ the dual co-frame in $T_x^*S^2$. Then the coframe   $\{\theta^1,\theta^2,\theta^3\}$ is easily obtained, the one-forms $\{\theta^1,\theta^2\}$ are the tautological forms of  $\alpha^1$, $\alpha^2$ on $T^*M$, while $\theta^3$ is the Levi-Civita connection form. 

Then this moving frame satisfies the structure equations
\begin{equation}
  \begin{split}
    & d\theta^1=\theta^2\wedge\theta^3\\
    & d\theta^2=\theta^3\wedge\theta^1\\
    & d\theta^3=G\theta^1\wedge\theta^2,
    \end{split}
  \end{equation}
  where $G$ is the Gauss curvature of $g$ regarded as function of one variable on $\Sigma$. Since we are using the odd function $h$ in \eqref{(1)}, the Gauss curvature is given by \eqref{G for g}. 
A quick glance at the formula \eqref{G for g} convinces that $g$ has positive Gauss curvature $G>0$ everywhere on $S^2$.
  
  Let us consider the following coframe changing
  \begin{equation}\label{Finsler coframe}
    \begin{pmatrix}
\omega^1\\ \omega^2 \\ \omega^3
\end{pmatrix}=
\begin{pmatrix}
  \sqrt{G} & 0 & 0 \\
  0 & 0 & -1\\
  0 & \sqrt{G} & 0
\end{pmatrix}
 \begin{pmatrix}
\theta^1\\ \theta^2 \\ \theta^3
\end{pmatrix}
\end{equation}
and by computing the structure equations of the coframe $\{\omega^1,\omega^2,\omega^3\}$, we obtain the structure equations
\begin{equation}
  \begin{split}
    & d\omega^1= -I \omega^1\wedge \omega^3+\omega^2\wedge\omega^3\\
    & d\omega^2= \omega^3\wedge\omega^1\\
     & d\omega^3=\omega^1\wedge\omega^2-J\omega^1\wedge\omega^3, 
  \end{split}
  \end{equation}
where we denote 
$$
I:=\frac{1}{2G\sqrt{G}}G_{\theta 2},\quad 
J:=-\frac{1}{2G\sqrt{G}}G_{\theta 1}.
$$
We point out that here, the function $G$ is the one in \eqref{G for g}.
  The subscripts here represent the directional derivatives with respect to the given co-frame, that is $df=f_{\theta 1}\theta^1+f_{\theta 2}\theta^2+f_{\theta 3}\theta^3$, for any smooth function $f:\Sigma\to \R$. 
  
  The 3-manifold $\Sigma$ endowed with the coframe $\{\omega^1,\omega^2,\omega^3
  \}$ is called a {\it generalized Finsler manifold} with the invariants $I$, $J$, $K=1$ (see \cite{4}, \cite{1}, etc.).
  
  Let us observe the following relations between the structures induced by $\{\theta^1,\theta^2,\theta^3\}$ and $\{\omega^1,\omega^2,\omega^3\}$:
  \begin{itemize}
  \item the geodesic foliation $\{\theta^1=0,\theta^3=0\}$ of $g$ coincides with the foliation $\{\omega^1=0,\omega^2=0\}$ (the indicatrix foliation);
  \item the indicatrix foliation $\{\theta^1=0,\theta^2=0\}$ of $g$ coincides with the foliation $\{\omega^1=0,\omega^3=0\}$ (the geodesic foliation).
  \end{itemize}

  Let us consider now the manifold of geodesics $M$ of the Zoll metric $(S^2,g)$. This is a 2-dimensional smooth manifold diffeomorphic to $S^2$. The following diagram is commutative
  
  \begin{center}
$\xymatrix
   { & \Sigma \ar[dl]_{\lambda} \ar[dr]^{\pi} \ar[r]^{\iota}&  TM& \\
     S^2  & & M
   }
$
  \end{center}
Here $\lambda:\Sigma\to S^2$ is the projection that identifies a leaf of the foliation $\{\theta^1=0,\theta^2=0\}$ with a point on $S^2$, and $\pi:\Sigma\to M$ is the projection that identifies a leaf of the foliation $\{\theta^1=0,\theta^3=0\}$ with a point on the manifold of geodesics $M$.

  It is clear that any Finsler structure on a manifold $M$ induces a generalized Finsler structure $\{\omega^1,\omega^2,\omega^3\}$ on the indicatrix bundle, where $\{\omega^1,\omega^2,\omega^3\}$ is the canonical orthonormal moving co-frame on $\Sigma$ determined by $F$
  (see \cite{4}, or any other textbook of Finsler geometry).

  However, clearly not every generalized Finsler structure leads to a Finsler metric. We will show however that the generalized Finsler structure constructed here actually gives a classical Finsler structure on $M$. 
  
  We recall here a fundamental result by Bryant (see \cite{3}).
\begin{proposition}\label{prop:Bryant}
	 The necessary and sufficient condition for a generalized Finsler surface  $
	 (\Sigma,\omega)$ to be realizable as a classical Finsler structure on a surface are
	 \begin{enumerate} 
	 	\item the leaves of the foliation $\{\omega^1=0,\ \omega^2=0\}$ are compact;
	 	\item it is amenable, i.e. the space of leaves of the foliation 
	 	$\{\omega^1=0,\ \omega^2=0\}$ is a differentiable manifold $M$;
	 	\item the canonical immersion $\iota:\Sigma\to TM$, given by 
	 	$\iota(u)=\pi_{*,u}(\hat{e}_2)$, is one-to-one on each $\pi$-fiber $\Sigma_x$,
	 \end{enumerate}
	 where we denote by $(\hat{e}_1,\hat{e}_2, \hat{e}_3)$ the dual frame of the coframing $(\omega^1,\omega^2,\omega^3)$.
\end{proposition}

In the same source it is pointed out that if for example the $\{\omega^1=0,\ \omega^2=0\}$ leaves are not 
compact, or even in the case they are, if they are ramified, or if the curves  $\Sigma_x$ winds around origin in $T_xM$, 
in any of these cases, the generalized Finsler surface structure is not realizable as a classical Finsler surface.

Obviously all conditions in Proposition \ref{prop:Bryant} are satisfied by $\Sigma$ endowed with the coframe \eqref{Finsler coframe} due to the properties of the Zoll metric $g$ (see \cite{3} and \cite{1}). It follows that the generalized Finsler metric defined by the coframe \eqref{Finsler coframe} is actually a classical Finsler metric on the manifold $M$ of geodesics of the Zoll metric constructed by using the function $h$ in \eqref{(1)}. 

\bigskip

Let us summarize our construction here.

\begin{itemize}
\item We start with a specific Zoll metric $(S^2,g)$ given by  \eqref{Zoll metric } with $h$ in \eqref{(1)}, with the Gauss curvature $G$ given by \eqref{G for g}. Observe that $G>0$.
  \item We consider the $g$-orthonormal coframe $\{\theta^1,\theta^2,\theta^3\}$ on the unit sphere bundle $\Sigma:=U^gS^2$.
\item We construct the coframe $\{\omega^1,\omega^2,\omega^3\}$ on $\Sigma$ by  \eqref{Finsler coframe}. Taking the exterior derivatives of the forms in this coframe, it is easy to see that this is a generalized Finsler structure on $\Sigma$ with the invariants $I$, $J$, $K=1$.
\item The generalized Finsler structure $(\Sigma;\omega^1,\omega^2,\omega^3)$ gives a classical Finsler structure on the manifold of geodesics $M$ of $g$, with flag curvature $K=1$.
  \item The Finsler structure obtained in this way is determined by the indicatrix curve in each tangent space $T_xM$, which is the curve \eqref{param eq indic ver1} obtained as the embedding of a $g$-geodesic in $TM$. 
\end{itemize}

This explaines why the Finsler metric constructed above is of constant flag curvature $K=1$.

\section{Appendix A.}

  Observe that in our case, $v_2$ in \eqref{param eq indic ver1} is given by
  \begin{equation*}
    v_2=\frac{1+\ve \sin^2 r\cos r}{\cos r}-\sqrt{\sin^2r-c^2}
    \int_R^r \frac{\sin s}{\cos^2 s}\frac{1}{\sqrt{\sin^2 s-c^2}}(1+2\ve \cos^3 s)ds.
    \end{equation*}

    The following formulas are useful for computing the above integral.
      \begin{equation*}
\int \frac{\sin s}{\cos^2 s}\frac{1}{\sqrt{\sin^2 s-c^2}}ds=\frac{1}{\cos^2 R}\frac{\sqrt{\sin^2 s-c^2}}{\cos s}+constant.
        \end{equation*}
        By integration by parts, we get
        \begin{equation*}
          \begin{split}
            & \int_R^r \frac{\sin s}{\cos^2 s}\frac{1}{\sqrt{\sin^2 s-c^2}}(1+2\ve \cos^3 s)ds\\
            &=\frac{1}{\cos^2R}
            \Bigg[
\frac{\sqrt{\sin^2r-c^2}}{\cos r}(1+2\ve \cos ^3r)+6\ve\int_R^r\sin s\cos s\sqrt{\sin^2s-c^2}\ ds
            \Bigg].
            \end{split}
          \end{equation*}
          \section{Appendix B.}

          The following formulas are useful when dealing with complex numbers.

          If $z=r(\cos\theta+i\sin\theta)$, then
          $$z^{\frac{1}{2}}=\pm\sqrt{r}\left(\cos\frac{\theta}{2}+i\sin \frac{\theta}{2}\right).$$

If $\bar{z}$ is the conjugate of $z$, then 
$$\bar{z}^{\frac{1}{2}}=\pm\sqrt{r}\left(\cos\frac{\theta}{2}-i\sin \frac{\theta}{2}\right).$$

It follows
$$
z^{\frac{1}{2}}+\bar{z}^{\frac{1}{2}}\in
\{\pm 2\sqrt{r}\cos\frac{\theta}{2},\pm 2i\sqrt{r}\sin\frac{\theta}{2}\}.
$$

Let us choose $z^{\frac{1}{2}}+\bar{z}^{\frac{1}{2}}=2\sqrt{r}\cos\frac{\theta}{2}$ in the previous formula. If $z=a+bi=re^{i\theta}$, $r:=\sqrt{a^2+b^2}$, $\tan \theta:=\frac{b}{a}$, then
$$
z^{\frac{1}{2}}+\bar{z}^{\frac{1}{2}}=\sqrt{2}\cdot \sqrt{a+\sqrt{a^2+b^2}}.
$$

Finally, if we let
$$a:=-\frac{1}{2}\left( P+Q-\frac{4\alpha}{3}\right),\quad b:=\frac{\sqrt{3}}{2}\left( P-Q
\right),
$$
then
\begin{equation*}
  \begin{split}
   & z^{\frac{1}{2}}+\bar{z}^{\frac{1}{2}}\\
   & =\sqrt{2}\cdot \sqrt{
     -\frac{1}{2}\left( P+Q-\frac{4\alpha}{3}\right)
     +\frac{1}{2}\sqrt{
\left( P+Q-\frac{4\alpha}{3}\right)^2+3(P-Q)^2
       }
     }.
\end{split}
\end{equation*}

Acknowledgements: The authors are greatful to K. Shibuya and U. Sombon for many useful discussions.

{\small {\sc Ioana Monica Ma\c sca}

{\sc Colegiul Na\c tional "Andrei \c Saguna", Bra\c sov}

E-mail: ioana.masca@yahoo.com}

\bigskip

{\small {\sc Sorin V. Sabau} (Corresponding author)

{\sc Departament of Biology, Tokai University, Sapporo, Japan}

E-mail: sorin@tokai.ac.jp}

\bigskip

{\small {\sc Hideo Shimada} 

{\sc Tokai University, Sapporo, Japan}

E-mail: shmd44@yahoo.co.jp}

\begin{thebibliography}{}

\bibitem[1]{2}
        D. {Bao}, S. S. {Chern}, Z. {Shen},
        {\it An introduction to Riemann - Finsler geometry}, Springer-Verlag, New York, 2000.
	
\bibitem[2]{4}
        D. {Bao}, C. {Robles}, Z. {Shen},
        {\it Zermelo navigation on Riemannian manifolds}, Journal of Differential Geometry, 66 (2004), 377--435.

\bibitem[3]{5}
        A. {Besse}, {\it Manifolds all of whose geodesics are closed}, Springer-Verlag, New York, 1978.		

\bibitem[4]{3}
        R. {Bryant}, {\it Some remarks on Finsler manifolds with constant flag curvature}, Houston Journal of Mathematics, Vol. 28, No. 2 (2002), 221-262.	
        	
\bibitem[5]{8} W. M.  Faucette, {\it A Geometric Interpretation of the Solution of the General Quartic Polynomial}, American Mathematical Monthly, 103 (1), (1996), 51--57.

\bibitem[6]{6}
  K.  {Kiyohara}, {\it Compact Liouville surfaces}, J. Math. Soc. Japan, Vol. 43, No. 3 (1991), 555--591.
  
       \bibitem[7]{1}
         K. Kiyohara, S. V. Sabau, K. Shibuya, {\it The geometry of a positive curved Zoll surface of revolution}, to appear in International Journal of Geometric Methods in Modern Physics
 (arXiv:1809.03138v2 [math.DG]).
         
\bibitem[8]{7}
        C. {Yu}, H. {Zhu}, {\it On a new class of Finsler metrics}, Differential Geometry and its Applications, 29 (2011), 244-254.
\end{thebibliography}
\end{document}